\documentclass[12pt]{article}
\usepackage{amssymb,amsmath,amscd,amsthm}

\newtheorem{theorem}{Theorem}[section]

\newtheorem{lemma}[theorem]{Lemma}

\date{}
\begin{document}
\title{ Intermittency for branching walks with heavy tails. }

\author{A. Getan \thanks{Department
of Mathematics and Statistics, University of North Carolina,
Charlotte, NC 28223, USA (agetan@uncc.edu).},  S. Molchanov \thanks{Department
of Mathematics and Statistics, University of North Carolina Charlotte, NC 28223, USA and National Research Univ., Higher School of Economics, Russion Federation. The work was partially supported  by the NSF grant DMS-1410547 (smolchan@uncc.edu).},  B. Vainberg \thanks{Department
of Mathematics and Statistics, University of North Carolina,
Charlotte, NC 28223, USA. The work was partially supported  by the NSF grant DMS-1410547 (brvainbe@uncc.edu). } \thanks{Corresponding author.}}

\maketitle

\begin{abstract}
Branching random walks on multidimensional lattice with heavy tails and a constant branching rate are considered. It is shown that under these conditions (heavy tails and constant rate), the front propagates exponentially fast, but the particles inside of the front are distributed very non-uniformly. The particles exhibit intermittent behavior in a large part of the region behind the front (i.e., the particles are concentrated only in very sparse spots there). The zone of non-intermittency (were particles are distributed relatively uniformly) extends with a power rate. This rate is found.
\end{abstract}



 \section{Introduction}

The  mathematical study of branching processes  goes back to the work of Galton  and
Watson \cite{wg} who were interested in the probabilities of long-term survival of family names. Later
similar mathematical models were used to describe the evolution of a
variety of biological populations, in genetics \cite{f1,f2,f3,hal}, and in the study of certain chemical
and nuclear reactions \cite{sem,haw}. The branching processes (in particular, branching diffusions)
play important role in the study of the evolution of various populations such as bacteria, cancer cells,
carriers of a particular gene, etc., where each member of the population may die or produce
offspring independently  of the rest.

In this paper, we describe the long-time behavior of a population in different regions of space
when, in addition to branching, the members of the population move in
space (random migration). In fact, we will consider discrete problems on random walks on the lattice $\mathbb Z^d$, not processes in $\mathbb R^d$, but our results can be extended to the latter case (with some changes).

The two main characteristics of the problem under consideration concern the rate of the branching and the probabilities of large jumps in $\mathbb Z^d$. There is an extensive literature on the branching random walks with compactly supported or fast decaying branching rates, see \cite{w,ah,br,aby1,aby2,by,ehk,mk, k,y2,or}. We assume that the branching rate $\nu$ is constant (the  total  population  grows exponentially  with
probability  one in this case). Problems of this type (with constant rate $\nu$ of splitting and zero mortality) originated from the KPP model, see the famous paper by Kolmogorov, Petrovski and Piskunov \cite{kpp}. The spreading of a new advanced gene (in $\mathbb R^d$) was studied in \cite{kpp}. We consider the problem on the lattice (not in $\mathbb R^d$), but there is a much more essential difference. We allow the random walk to have long jumps. To be more exact, we consider random walks with heavy tails. This combination of heavy tails with the constant branching rate creates the following effect: the front of the population propagates exponentially fast, but the number of particles inside of the front is distributed very  non-uniformly. The latter property is referred to as the intermittency for the number of particles, and is the main object of our study.

Let us formulate the problem more precisely.  Branching walks on the lattice $\mathbb Z^d$ are considered. We assume that a particle located at $y \in\mathbb Z^d$ at a given time $t$ remains at the same point for a random exponentially distributed time (with parameter $1$), and then jumps to a new position  $y+z$ with probability $a(z)$. In addition, each particle located at $y\in \mathbb Z^d$  splits with the rate $\nu dt$ into two particles located at the same point $y$. All the particles behave independently of each other and according to the same law.

Obviously,
\begin{equation}\label{00}
a(z)\geq 0, \quad \sum_{z\in \mathbb Z^d}a(z)=1.
\end{equation}
We also assume that the distribution of the jumps is symmetric, i.e., $a(z)=a(-z)$, which implies that  $L=L^*$ in $l^2 (\mathbb Z^d)$. We assume that $a$ has the following behavior at infinity
\begin{equation}\label{ainf}
a(z)= \frac{a_0(\dot{z})}{|z|^{d+\alpha}}(1+o(1)), \quad |z|\to \infty , \quad \dot{z}= \frac{z}{|z|},
\end{equation}
 {\rm where}
   \begin{equation}\label{ainf1}
a_0(\dot{z})> \delta > 0, \quad 0<\alpha <2.
  \end{equation}
In fact, later we will specify the asymptotics of $a(z)$ in more detail (see (\ref{005})).

Many essential characteristics of the random walk (or process in $R^d$) depend on the rate of decay of the distribution at infinity. Conditions (\ref{ainf}), (\ref{ainf1}) correspond to a walk with a heavy tail. The second moments do not exist in this case. The decay (\ref{ainf}) with $\alpha>2$ defines walks with a moderate tail. A process has a light tail if the distribution decays at infinity so fast that its Fourier transform $\widehat{a}(k)$ is analytic in $k$ when $|{\rm Im}k|<\delta$ with some $\delta >0$. The second moments are well defined in the latter two cases. Thus the walk under consideration admits very long jumps whose probabilities are not very small.

Let $n(t,x,y)$ be the  number of particles at a point $x \in Z^d$ at the moment $t \geq 0$ under condition that the process starts at initial moment $t=0$ with a single particle located at a point $y \in\mathbb Z^d$,
i.e., $n(0,x,y)=\delta(x-y)$.

Without duplication, the initial particle would perform the symmetric random walk $X(t)$ with the generator
\[
(\mathcal L_y f)(y)= \sum_{z\in \mathbb Z^d}[f(y + z)-f(y)]a(z),
\]
which acts on the space $\textit{l}^2(\mathbb Z^d).$

If $\nu>0$ and $m_1= m_1(t,x,y)$ is the expected value (the first moment) of the random variable $n(t,x,y)$, then $m_1$ satisfies the relations:
\begin{equation}\label{n00m1}
\frac{\partial{m_1}}{\partial t} (t,x,y)=  \mathcal L_y m_1(t,x,y)+ \nu~  m_1(t,x,y), \quad   t\geq 0.
\end{equation}
\begin{equation*}
     m_1(0,x,y)=\delta(x-y).\quad
\end{equation*}

In order to derive (\ref{n00m1}), we evaluate $ m_1(t+\Delta t,x,y) $ by splitting the time interval $ (0, t + \Delta t) $ into two successive intervals of lengths $\Delta t$ and $ t $. Then

\begin{equation} \label{000}
m_1(t+\Delta t,x,y)\sim \sum_{z\in \mathbb Z^d}{a(z)\Delta t m_1(t,x,y+z)}+2 \nu \Delta t m_1(t,x,y)
\end{equation}
\[~~~~~~~~~~~+(1-\Delta t -\nu\Delta t)m_1(t,x,y).~~~~~~~~~~~\]

The first term on the right side of $ (\ref{000} )$   is the sum of the probabilities to jump from $ y $ to $ y+z $ during the time $\Delta t$ (these probabilities are $ a(z)\Delta t $) multiplied by the expectation $ m_1(t,x,y+z)$ for the number of particles at $x$ when the walk starts at a single point $ y+z \in \mathbb Z^d$.
The second term describes the probability $\nu \Delta t $ of branching during the time $ \Delta t $ multiplied by the expected number of particles at $ x $ that are descendants of both the original and the new particles at $ y $, which is $ 2m_1(t,x,y)$.
The last term is the contribution to the expectation of the number of particles at $y$ from the event that the particle stays at $ x $ without branching and jumping during the time interval $\Delta t $.

We subtract $m_1(t,x,y)$ from both sides of equation $ (\ref{000} )$ above, divide by $\Delta t$, and pass to the limit $\Delta t \rightarrow 0 $. This implies (\ref{n00m1}).

Let us derive the equation for the second moment,
$m_2 (t,x,y)= E(n^2(t,x,y)) $.
We again consider the time interval $(0,t+ \Delta t) $ and split it in two successive intervals of lengths $ \Delta t $ and $t $. Then
\begin{equation*}
m_2(t+\Delta t,x,y)\sim \sum_{z\in \mathbb Z^d}{a(z)\Delta t~ m_2(t,x,y+z)+ \nu~ \Delta t~ E(n_1 + n_2)^2}
\end{equation*}
\begin{equation}\label{00m2}
~~~~~~~~~~~~~~~~~+(1-\Delta t -\nu~\Delta t)~m_2(t,x,y).
\end{equation}
Here the terms on the right of (\ref{00m2})  are similar to the terms in (\ref{000}),
$n_1(t,x,y)$ is the number of particles in $ x $ that are descendants of the original particle, and $n_2(t,x,y)$ is the number of particles in $ x$ that are descendants of the newly born particle.
We use the fact that
$E(n_1 + n_2)^2=E(n_1 )^2 + E( n_2)^2 + 2E(n_1) E(n_2)
= 2m_2(t,x,y)+2m_1^2(t,x,y)$.
Then we subtract $ m_2(t,x,y) $ from both sides of (\ref{00m2} ), divide by $ \Delta t $ and pass to the limit as $ \Delta t \ \rightarrow 0 $. This implies:

\begin{equation}\label{01m2}
\frac{\partial{m_2}}{\partial t}(t,x,y)= (\mathcal L_y +\nu)m_2(t,x,y)+ 2 \nu~ m_1^2 (t,x,y),
\end{equation}
\begin{equation*}
  m_2 (0,x,y)= \delta(x-y).
\end{equation*}

Denote by $m_1(t,x),~m_2(t,x)$ the first and second moments for the number of particles at the point $x$ when the process starts at a single particle at the origin:
$$
m_1(t,x):=m_1(t,x,0),~~m_2(t,x):=m_2(t,x,0)
$$
Obviously, functions $m_1(t,x,y),m_2(t,x,y)$ depend on $x-y$, not on $x$ and $y$ separately, i.e.,
$$
m_1(t,x,y)=m_1(t,x-y,0),~~m_2(t,x,y)=m_2(t,x-y,0),
$$
and therefore, $m_1(t,x,y),m_2(t,x,y)$ are known as soon as  $m_1(t,x),m_2(t,x)$ are found. Note also that the symmetry of $a(z)$ implies that $\mathcal L_yf
(x-y)=\mathcal Lf(x-y)$, where
$$
\mathcal Lf=\mathcal L_xf=\sum_{z\in \mathbb Z^d}[f(x + z)-f(x)]a(z)
$$
Hence $m_1(t,x),~m_2(t,x)$ satisfy the relations
\begin{equation}\label{n00m1a}
\frac{\partial{m_1}}{\partial t} (t,x)=  (\mathcal L+\nu) m_1(t,x), \quad   t\geq 0; \quad m_1(0,x)=\delta(x).
\end{equation}
\begin{equation}\label{01m2a}
\frac{\partial{m_2}}{\partial t}(t,x)= (\mathcal L+\nu)m_2(t,x)+ 2 \nu~ m_1^2 (t,x), \quad   t\geq 0;\quad m_2(0,x)=\delta(x).
\end{equation}

The region in $\mathbb Z^d$ that separates the large and small values of $m_1(t,x)$ is called {\it the front}. To be more exact, we define the front $F=F(t)$ as the boundary of the set $\{x: m_1(t,x)\geq 1\}$. The boundary consists of points $x$ that have neighbors $x+e_{\pm },~|e_{\pm }|=1$ or $e_{\pm }=0$, such that  $m_1(t,x+e_+)\geq 1, ~m_1(t,x+e_-)< 1$.

The notion of {\it intermittency} (or intermittent random fields) is popular in natural  sciences (astrophysics, biology, etc). From the qualitative point of view, intermittent random fields are distinguished by the formation of sparse spatial structures such as high peaks, clumps, patches, etc.,
giving the main contribution to the process in the medium. For instance, the magnetic field of the Sun is highly intermittent as almost all its energy is concentrated in the black spots, which cover only a very small part of the surface of the Sun. Many bio-populations also exhibit strong clumping (clustering).

Intermittency is a well developed non-uniformity. For physicists, the magnetic field of the Sun is intermittent since, say, $99\%$ of its magnetic energy is concentrated on less than $1\%$ of the surface. For mathematicians, $0.1, 0.01$ or $10^{-6}$ are not necessarily small numbers, and a limiting process must be considered instead. The definition of intermittency based on the progressive growth of the statistical moments was proposed in the review \cite{zel}, a more formal presentation can be found in \cite{SantFl}. In the the simplest form, a field $n(t,x),~x\in \mathbb Z^d, $ is intermittent as $t\to\infty$ on a non-decreasing family of sets $D(t)$ if
\[
\lim_{t\to\infty}\frac{En^2(t,x)}{(En(t,x))^2}=\infty
\]
uniformly in $x\in D(t)$.

Let us illustrate this definition with the following example. Let $n(t,x),~x\in \mathbb Z^d, $ be independent identically distributed r. v. and
\[
P\{n(t,x)=0\}=1-\frac{p_0}{t}, ~~P\{n(t,x)=t\}=\frac{p_0}{t}, ~~t\geq p_0.
\]
Then $En(t,x)=p_0$, i.e., the density of the population is constant in time, and
\[
p_0=\lim_{L\to\infty}\frac{\sum_{|x|\leq L}n(t,x)}{(2L)^d} \quad ({\rm the ~law ~ of ~ the ~ large~  numbers}).
\]
However,
\[
\frac{En^2(t,x)}{(En(t,x))^2}=p_0t\to\infty,
\]
and similarly for $\frac{En^m(t,x)}{(En(t,x))^m}, ~m\geq 2$. Thus, the family of the fields $n(t,x)$ is intermittent as $t\to\infty$. It is clear that the population $n(t,x)$ is supported for large $t$ on a subset of $\mathbb Z^d$ with relative volume $p_0/t\to 0$ as $t\to\infty$. The independence of $n(t,x)$ is not important here and can be replaced by some kind of weak dependence (ergodicity), see \cite{SantFl}.


The main two results of the present paper concern the propagation of the front of the branching random walk that starts at the origin, and the intermittency of the distribution of the particles on and behind the front. These results are proved under the following assumption that is a little more restrictive than (\ref{ainf}).
Namely, we assume that:

\begin{equation}\label{005}
a(z)= \sum_{j=0}^{d+\epsilon}\frac{a_j(\dot{z})}{|z|^{d+\alpha+j}}+O(\frac{1}{|z|^{2d+\alpha+1+\epsilon}}), \quad |z|\to\infty, ~~~ \alpha\in (0,2),
\end{equation}
where $\dot{z}=z/|z|,$
\begin{equation*}
a_j\in C^{d+1-j+\epsilon}(S^{d-1}), \quad a_0(\dot{z})>\delta>0,
\end{equation*}
and $\epsilon=1$ if $\alpha=1$, $\epsilon=0$ otherwise.
\begin{theorem}\label{t1}
Let (\ref{005}) hold. Then the following asymptotics holds for the points $x\in F(t)$ on the front $F(t)$:
\begin{equation} \label{f1}
|x|=[a_0(\dot{x})t]^\frac{1}{d+\alpha}e^{\frac{\nu}{d+\alpha}t}(1+O(t^{-1})),  \quad  t\to\infty.
\end{equation}
\end{theorem}
This theorem is an immediate consequence of the global limit theorem proved in \cite{arxive} (see Theorem \ref{mt} below and the remark after it).
While the front propagates exponentially fast, the particles are distributed very non-uniformly on the front and at any exponential in time distance from the original particle. In fact, we will find the exact boundary for non-intermittency, and this boundary propagates with a power rate. Let
\begin{equation} \label{gamma}
\gamma=\frac{2\alpha +d}{\alpha(\alpha+d)}.
\end{equation}
The following statement will be proved in the next section.
\begin{theorem}\label{t2}
Let (\ref{005}) hold. Then

1) The ratio $\frac{m_2 (t,x)}{m_1^2 (t,x)}$ is uniformly bounded in each ball $|x|<Bt^\gamma$ when $t\to\infty$, i.e., the random variable $n$ is non-intermittent there.

2) For each domain $\Omega_\varepsilon(t)=\{x:|x|>t^{\gamma+\varepsilon}\},~\varepsilon>0$, we have $\frac{m_2 (t,x)}{m_1^2 (t,x)}\to\infty$ uniformly in $x\in \Omega_\varepsilon(t)$ as $t\to\infty$, i.e., $n$ is intermittent in $\Omega_\varepsilon(t)$.
\end{theorem}

\section{Proof of the main results.}
The function $m_1(t,x)$ with $\nu = 0$  will be denoted by $p = p(t,x)$. Then
\begin{equation}\label{00p}
\frac{\partial p}{\partial t}(t,x)=\mathcal L p(t,x)~, ~~t\geq 0~; \quad p(0,x)=\delta(x),
\end{equation}
and
\begin{equation}\label{p1m}
  m_1(t,x)= e^{\nu t} p(t,x).
\end{equation}

Consider the homogeneous (of order $-d-\alpha$) distribution in $\mathbb R^d$ that is equal to $a_0(\dot{x})|x|^{-d-\alpha}$ when $x\neq 0$ (compare with the first term of asymptotics (\ref{ainf})). It was shown in \cite{arxive} that the Fourier transform of this distribution is
the homogeneous function   $-b_0(\dot{\sigma})|\sigma|^{\alpha}$ in $R^d$, where $\sigma$ is the dual variable to $x,~\dot{\sigma}=\frac{\sigma}{|\sigma|},$ and
\begin{equation}\label{007}
b_0(\dot{\sigma})=-\Gamma(-\alpha)\cos\frac{\alpha\pi}{2}\int_{S^{d-1}}a_0(\dot{x})|(\dot{x},\dot{\sigma})|^{\alpha}dS_{\dot{x}}>0.
\end{equation}
Here $\Gamma$ is the gamma-function.

The following global limit theorem for random walks with heavy tails obtained in \cite{arxive} is a key point in the proof of the results stated in the introduction.

\begin{theorem} \label{mt} Let (\ref{ainf}) hold. Then

 (i) the following asymptotis holds for $p$
\begin{equation}\label{009}
 p(t,x)=\frac{1}{t^{d/\alpha}}S(\frac{x}{t^{1/\alpha}})(1+o(1)) \quad {\it when} \quad x\in \mathbb Z^d, \quad |x|+t\to\infty,
\end{equation}
were
$S(y)=\frac{1}{(2\pi)^{d}}\int_{R^d} e^{i(\sigma,y)-b_0(\dot{\sigma})|\sigma|^{\alpha}}d\sigma>0$ is the stable density $S=S_{\alpha,a_0}(y)$, which depends on $\alpha\in(0,2)$ and $a_0$, and $b_0$ is defined in (\ref{007}).

(ii)\quad If $\frac{|x|}{t^{1/\alpha}}\to\infty, ~|x|\geq 1,$ then the previous statement can be specified as follows:

\begin{equation}\label{0010}
p(t,x)=\frac{a_0(\dot{x})}{t^{d/\alpha}}(\frac{t^{1/\alpha}}{|x|})^{d+\alpha}(1+o(1))=\frac{a_0(\dot{x})t}{|x|^{d+\alpha}}(1+o(1)).
 \end{equation}
\end{theorem}
\noindent{\bf Remarks.} 1. There is a misprint in the statement of the theorem in \cite{arxive}: $S(\frac{|x|}{t^{1/\alpha}})$ must be replaced by $S(\frac{x}{t^{1/\alpha}})$.

2. This theorem and (\ref{p1m}) imply Theorem \ref{t1}. The remaining part of the section will be devoted to the proof of Theorem \ref{t2}.

For a given function $f=f(x),~x\in \mathbb Z^d$, denote by $\widehat{f}(\sigma)$ the periodic in $\sigma$ function that is the Fourier series with coefficients $f(x), x\in \mathbb Z^d$, i.e.,
\[
\widehat{f}(\sigma)=\sum_{x\in \mathbb Z^d}f(x)e^{-i(x,\sigma)}, \quad f(x)=\frac{1}{(2\pi)^d}\int_{[-\pi,\pi]^d}\widehat{f}(\sigma)e^{i(x,\sigma)}d\sigma.
\]
In particular, from (\ref{00p}) it follows that
 \begin{equation}\label{003}
 p(t,x)= \frac{1}{(2\pi)^d}\int_{[-\pi,\pi]^d} {\widehat{p}(t,\sigma)e^{i(\sigma,x)}d\sigma}= \frac{1}{(2\pi)^d}\int_{[-\pi,\pi]^d}{e^{[ \widehat{a}(\sigma)-1]t + i(\sigma,x)}d \sigma}.
 \end{equation}
 The following properties of $\widehat{a}(\sigma)$ follow immediately from properties of $a(x)$:
\begin{equation}\label{008}
\widehat{a}(-\sigma)=\widehat{a}(\sigma);\quad  ~-1<\widehat{a}(\sigma)<1, ~~ 0\neq \sigma\in T^d.
\end{equation}
The second part in (\ref{008}) follows from (\ref{00}) provided that for each $\sigma\in T^d, ~\sigma\neq 0,$ there is a point $z\in Z^d$ where $e^{-i(z,\sigma)}\neq 1$ and $a(z)\neq 0$. Such points $z$ exist due to (\ref{005}).

We will need the following lemma.
\begin{lemma}\label{ppol}
 Function $p(t,x)$ is strictly positive for all $x \in \mathbb Z^d$, $t>0$.
\end{lemma}
{\bf Proof.}
Denote by $a_n (x)$ the convolution of $n$ copies of $a(x)$:
\begin{equation}\label{0*1}
 a_n (x): = a(x)*a(x)* ... *a(x),
\end{equation}
where  $ a(x)*b(x)= \sum_{z\in \mathbb Z^d}{a(x-z)*b(z)} $.

We multiply both sides of (\ref{0*1}) by $e^{(-i\sigma x )}$ and take the sum in $x \in \mathbb Z^d$. This implies
$\widehat{a}_n(\sigma) =[\widehat{a}(\sigma)]^n.$  From here, (\ref{00}), and (\ref{008}) it follows that
 $|\widehat{a}_n(\sigma)|\leq 1$. This allows us to write $p(t,x)$ as follows:
\begin{equation*}
p(t,x)= \frac{1}{(2\pi)^d}\int_{[-\pi,\pi]^d}{e^{[\widehat{a}(\sigma)-1]t + i(\sigma,x)}d \sigma}= \frac{e^{-t}}{(2\pi)^d}\int_{[-\pi,\pi]^d}{[1+\sum_{n=1}^{\infty}{\frac{[\widehat{a}(\sigma)]^n}{n!} }t^n ]~ e^{i(\sigma,x)} d\sigma}.
\end{equation*}

Thus

\begin{equation}\label{0*3}
p(t,x) = e^{-t}[\delta(x) +\sum_{n=1}^{\infty} \frac{a_n(x)}{n!}t^n  ].
 \end{equation}

 Since $a(x)\geq0$, all the convolutions $a_n (x)$ are non-negative. Hence (\ref{0*3}) will imply the statement of the lemma if we show that $a_2 (x)$ is strictly positive for all $x \in \mathbb Z^d$.

 We have:
 \begin{equation}\label{0*4}
 a_2(x)=\sum_{z \in \mathbb Z^d}{a(x-z)a(z)}.
 \end{equation}

 Here $a\geq0$, and from (\ref{005}) it follows that the terms in (\ref{0*4}) are positive for each fixed $x$ if $z$ is large enough. Thus $a_2 (x)>0$, and the proof of the lemma is complete.

 \qed

 Let us formulate a simplified version of Theorem \ref{mt} that will be combined with (\ref{p1m}) and will be easer to use than Theorem \ref{mt}.

\noindent{\it Definition.} Functions $a$ and $b$ will be called equivalent, and it will be denoted by $a\asymp b$, if there exist two constants $c_1$ and $c_2$ such that,  $c_1 b<a<c_2 b$.

\begin{lemma} \label{linf}
Let (\ref{ainf}) hold. Then for arbitrary $a_2\geq a_1>0$, the following relations hold
\begin{equation}\label{t2i}
 (i)~~~|m_1(t,x)|\asymp  \frac{t}{|x|^{d+\alpha}}e^{\nu t} ~~   {\it  when} ~~|x|> a_1t^{\frac{1}{\alpha}},
 \end{equation}
\begin{equation}\label{t2ii}
 (ii)~~~|m_1(t,x)|\asymp
t^{-\frac{d}{\alpha}}e^{\nu t} ~~{\it when}~~ |x|\leq a_2t^{\frac{1}{\alpha}},~~t>\varepsilon>0,
 \end{equation}
\begin{equation}\label{t2iii}
 (iii)~~~|m_1(t,x)|\asymp 1 ~~{\it when} ~~x=0,~ t<1.
 \end{equation}
 \end{lemma}
\noindent{\bf Remarks.} 1) One can use any of the estimates (\ref{t2i}), (\ref{t2ii}) in the intermediate zone $a_1t^{\frac{1}{\alpha}}<|x|\leq a_2t^{\frac{1}{\alpha}},~~t>\varepsilon>0.$ The right-hand sides of these estimates are equivalent in this intermediate zone.

2) From the last two relations it follows that the following estimate holds for all $t>0$:
\begin{equation}\label{t20}
 |m_1(t,x)|\asymp (t+1)^{-\frac{d}{\alpha}}e^{\nu t} ~~{\rm when}~~ |x|\leq a_2t^{\frac{1}{\alpha}}.
 \end{equation}
{\bf Proof.} Due to (\ref{p1m}), it is enough to prove the lemma when $\nu=0$ and $m_1$ is replaced by $p$.

The last statement of the lemma follows immediately from (\ref{003}).
In order to prove the first two statements, we split the region $x \in \mathbb Z^d, t \geq 0, $ into three subregions,
$U_1,U_2,$ and $ U_3,$ where $U_1$  is the region defined by the inequality $  |x|>At^\frac{1}{\alpha} $ with $A$ so large that the remainder term in  (\ref{0010}) is less than $1/2.$
Note that $|x|>At^{\frac{1}{\alpha}}$ implies that $x\neq 0.$ Thus $|x|\geq 1$ in $U_1 $ since $x$ is a point on the lattice. Hence (\ref{0010}) implies (\ref{t2i}) with $a_1=A$.
Let $U_2$ be defined by the inequalities $ |x|\leq At^{\frac{1}{\alpha}},~ |x|+t\geq B   $, where $B$ is chosen so large that the remainder term in  $(\ref{009})$ is less than $1/2$  when $|x|+t>B.$
Since  $\frac{|x|}{t^ {\frac{1}{\alpha}}}\leq A $ is bounded in $U_2 $ and function $S(y)$ is positive and continuous, it follows that $S(\frac{x}{t^ {\frac{1}{\alpha}}})$ in $(\ref{009})$
 has upper and lower positive bounds in $U_2$. Thus $(\ref{009})$ implies that  $ p(t,x) \asymp t^{-\frac{d}{\alpha}}$ in $U_2$.

Now consider the region
\[
U_3= \begin{Large} \{ |x|\leq At^{\frac{1}{\alpha}}  ,~ |x|+t\leq B \} \backslash \{(t,x): ~x=0,~t<1 \} \end{Large}.
\]
Region $U_3$ is bounded, and $t\geq\delta>0$ there. Since $p$ is continuous, from Lemma \ref{ppol} it follows that $p$
has positive lower and upper bounds on $U_3$. Similar bounds are valid for the function $
 t^{-\frac{d}{\alpha}}$.
Thus $ p(t,x)  \asymp t^{-\frac{d}{\alpha}} $  on $U_2 \bigcup U_3$, and the second statement of Lemma \ref{linf} is proved with $a_2=A$. Estimates (\ref{t2i}), (\ref{t2ii}) were proved  with $a_1=a_2=A$. Their validity with arbitrary $a_1,~a_2$ follows from the equivalency relation stated in the remark above.

\qed

\noindent {\bf Proof of Theorem \ref{t2}.}
Since $m_1$ is the Green function for the operator $\mathcal L+\nu$, the Duhamel principle implies that the solution $m_2$ of (\ref{0014}) has the form
\begin{equation}\label{0014}
 m_2 (t,x)= m_1 (t,x)+2\nu\int_{0}^{t}\sum_{z\in \mathbb Z^d}m_1 (t-s,x-z)  m_1^2 (s,z) ds.
\end{equation}

We will start the proof of the theorem with the second statement (about the intermittency). Since $m_1,m_2\geq 0$, we have
\begin{equation}\label{0014a}
 m_2 (t,x)\geq 2\nu\int_{t-1}^{t}\sum_{z\in \mathbb Z^d:~|z|\leq \frac{1}{2}t^{\frac{1}{\alpha}}}m_1 (t-s,x-z)  m_1^2 (s,z) ds,
\end{equation}
where $x\in \Omega_\varepsilon(t), ~ t\to\infty.$ We apply estimate (\ref{t2i}) to the first factor under the summation sign and the estimate  (\ref{t2ii}) to the second factor. Taking into account that $t-s\leq 1$ and $|x-z|\asymp |x|$ in (\ref{0014a}), we obtain that
\begin{equation*}
 m _2 (t,x)\geq  C \int_{t-1}^{t}
\sum_{z\in \mathbb Z^d:~|z|\leq \frac{1}{2}t^{\frac{1}{\alpha}}}\frac{ 1}{|x |^{d+\alpha~}} ~s^{-\frac{2d}{\alpha}}e^{2\nu t} ds\geq
C \int_{t-1}^{t}
\frac{t^{\frac{d}{\alpha}}}{|x |^{d+\alpha~}} ~s^{-\frac{2d}{\alpha}}e^{2\nu t} ds.
\end{equation*}
Hence
 \begin{equation*}
  m_2 (t,x)\geq   \frac{Ce^{2\nu t}}{t^{\frac{d}{\alpha}}|x|^{d+\alpha}}, \quad x\in \Omega_\varepsilon(t), ~~~ t\to\infty.
\end{equation*}
On the other hand, from  (\ref{t2i}) it follows that
\begin{equation*}
m_1^2 (t,x)  \leq   \frac{C t^2}{|x|^{2(d+\alpha)}}  e^{2\nu t}, \quad x\in \Omega_\varepsilon(t), ~~~ t\to\infty.
\end{equation*}
Hence, if $x\in \Omega_\varepsilon(t), ~ t\to\infty$, then
\begin{equation*}  \frac{m_2 (t,x)}{m_1^2 (t,x)} \geq C \frac{|x|^{ (d+\alpha)}}{t^{2+d/\alpha}}  \geq C\frac{t^{(\gamma+\varepsilon)(d+\alpha )}}{t^{\frac{2\alpha+d}{\alpha}}}=Ct^{\varepsilon(d+\alpha)}\to\infty \quad {\rm as} \quad t\to\infty.
\end{equation*}
The second statement of the theorem is proved.

Let us prove the first statement. Note that the ball $|x|\leq Bt^\gamma$ is located behind the front, and Lemma \ref{linf} implies that $m_1\to\infty$ in this ball as $t\to\infty$. Hence $(\ref{0014})$ implies that the first statement of the theorem will be proved (and the proof of the theorem will be complete) as soon as we show that

\begin{equation}\label{0019}
I:=\frac{1}{m_1^2(t,x)}
\int_{0}^{t}
\sum_{z\in \mathbb Z^d}m_1 (t-s,x-z) m_1^2 (s,z)~ds<C,  \quad ~  |x|\leq Bt^\gamma, \quad t\to\infty.
  \end{equation}

  In order to estimate the left hand-side of (\ref{0019}), we split $\mathbb Z^d$ in $(\ref{0019})$ into four sets separated by the two spheres:
 $|z|= s^{\frac{1}{\alpha}}$ and $|x-z|=(t-s)^{\frac{1}{\alpha}}$. Let $P_1=P_1(s), P_2=P_2(x,t-s)$ be (bounded) sets of points $z\in \mathbb Z^d$ located inside or at the boundary of the first (respectively, second) sphere defined above, i.e.,
\[
P_1=\{z\in \mathbb Z^d: |z|\leq s^{\frac{1}{\alpha}}\}, \quad P_2=\{z\in \mathbb Z^d: |x-z|\leq(t-s)^{\frac{1}{\alpha}}\}
\]
Denote by $D_{i,j}=D_{i,j}(s,t,x)$ the following sets of points $z\in \mathbb Z^d$ with $0\leq s\leq t,~x\in\mathbb Z^d$:

$D_{11}$ is the the set of points $z\in \mathbb Z^d$ located inside or at the boundary of both spheres, i.e., $D_{11}=P_1\bigcap P_2.$

$D_{22}$ is the the set of points $z\in \mathbb Z^d$ located outside of both spheres, i.e., $D_{22}=\mathbb Z^d\setminus ( P_1\bigcup P_2).$

$D_{12}$ is the set of points $z\in \mathbb Z^d$ located outside of the first sphere, but inside of the second one or on its boundary, i.e., $D_{12}=(\mathbb Z^d\setminus P_1)\bigcup P_2.$

$D_{21}$ is the the set of points $z\in \mathbb Z^d$ located inside of the first  sphere or on its boundary, but outside of the second one, i.e., $D_{21}=(\mathbb Z^d\setminus P_2)\bigcup P_1.$

Respectively, $I$ can be written as $I=\sum_{i,j=1}^2I_{ij}$, where $I_{ij}$ is defined by (\ref{0019}) with the summation extended oved $D_{ij}$ instead of $\mathbb Z^d$. We are going to estimate each of the terms $I_{ij}$ when $|x|\leq Bt^\gamma, ~ t\to\infty$.

{\it 1) Estimate on $I_{12}$}. We will obtain this estimate separately for $2t^{\frac{1}{\alpha}}<|x|\leq Bt^\gamma$ and for $|x|\leq 2t^{\frac{1}{\alpha}}$.
Note that the two balls $P_1$ and $P_2$ are separated in the first case (when $0\leq s\leq t$) and they may intersect each other in the second case. Consider the first case of $|x|>2t^{\frac{1}{\alpha}}$. Then $|m_1^2 (t,x)| \asymp  \frac{t^2}{|x|^{2(d+\alpha)}} e^{2t\nu}$ (due to (\ref{t2i})) and relations (\ref{t20}), (\ref{t2i}) hold for the first and second factors under the summation sign in (\ref{0019}), respectively. Thus

\begin{equation}\label{0020}
I_{12} \leq  \frac{C|x|^{2( d+\alpha)}}{t^2 e^{2t\nu} }
  \int _{0}^{t}{
  \sum
 _{z\in D_{12} }{\frac{s^2 e^{\nu(t+s)}ds}{~[(t-s)+1]^{\frac{d}{\alpha}}|z|^{2(d+\alpha)}}}~}, \quad 2t ^{\frac{1}{\alpha}}<|x|\leq Bt^\gamma, ~ ~ t\to\infty.
\end{equation}

 Since $|x-z|\leq (t-s) ^ {\frac{1}{\alpha}} \leq t ^{\frac{1}{\alpha}} $ in $D_{12}$, inequality $|x|>2t ^{\frac{1}{\alpha}}$ implies that $|x-z|\leq \frac{1}{2} |x|$, and therefore $|z| \geq \frac{1}{2}|x|$ in (\ref{0020}).
 Hence
we can replace $z$ by $x$ there. After that, summation in (\ref{0020}) is applied to $z$-independent terms. Hence, the summation sign can be replaced by a factor $\kappa$ that estimates the number of terms in the sum from above. Obviously, $\kappa\leq C(A^{d}+1)$, where $A$ is the radius of the ball $P_2$, i.e., $ \kappa\leq C[(t-s)^{\frac{d}{\alpha}}+1],$ and
\begin{equation}\label{D12C1}
I_{12} \leq \frac{C}{t^2 e^{2\nu t}}\int_{0}^{t} \frac{ [(t-s)^{\frac{d}{\alpha}}+1]s^2 e^{\nu(t+s)} ds}{(t-s+1)^{\frac{d}{\alpha}}} \leq  \frac{C_1}{t^2 e^{ \nu t}}\int_{0}^{t} s^2 e^{\nu s} ds \leq C_2< \infty
\end{equation}
when $2t ^{\frac{1}{\alpha}}<|x|\leq Bt^\gamma, ~ ~ t\to\infty.$

  Now let $|x|\leq 2t^{\frac{1}{\alpha}}$. Then
    $|m_1^2 (t,x)| \asymp
 \frac{1}{t^{\frac{2d}{\alpha}}} e^{2t\nu}$ for   $t\geq 1$ (due to (\ref{t2ii})), and  (\ref{0020}) must be replaced by
 \begin{equation}\label{2D12}
  I_{12} \leq  \frac{C t^{\frac{2d}{\alpha}}}{ e^{2\nu t}}
  \int _{0}^{t}{
  \sum_{z\in D_{12} }{\frac{s^2 e^{\nu(t+s)}ds}{(t-s+1)^{\frac{d}{\alpha}}|z|^{2(d+\alpha)}}}~}, \quad |x|\leq 2t ^{\frac{1}{\alpha}}, ~ ~ t\to\infty.
 \end{equation}.

We split the right-hand side above into two parts $I'+I''$  by writing the interval $(0,t)$ as the union of $(0,t/2)$ and $(t/2,t)$. Since $D_{12}$ does not contain points of the first ball $P_1$, $z\neq 0$ there. We have $|z|\geq 1$ for all other points of the lattice. Thus we can replace $|z|$ in $I'$ by $1$ and replace summation by factor $\kappa$. This leads to the estimate
$$
I'\leq\frac{C t^{\frac{2d}{\alpha}}}{ e^{2\nu t}}
  \int _{0}^{t/2} \frac{ [(t-s)^{\frac{d}{\alpha}}+1]s^2 e^{\nu(t+s)} ds}{(t-s+1)^{\frac{d}{\alpha}}} \leq \frac{C t^{\frac{2d}{\alpha}}}{ e^{\nu t}}
  \int _{0}^{t/2}s^2 e^{\nu s} ds,
  $$
 where the right hand side decays exponentially as $t\to \infty$. Thus $I'\leq C< \infty.$

 Now consider $ I'' $. In this case, $t/2<s<t$ implies that $ |z|\geq s^{\frac{1}{\alpha}}\geq {(t/2)}^{\frac{1}{\alpha}} $, and we can replace $|z|$ by  ${(t/2)}^{\frac{1}{\alpha}}$ in (\ref{2D12}).
 Then we replace the summation in z by the factor $\kappa$  as before.
 This leads to the following estimate:
 \begin{equation*}
 I''\leq\frac{C t^{\frac{2d}{\alpha}}}{ e^{2\nu t}t^{2(d+\alpha)/\alpha}}
  \int_{t/2}^{t}s^2 e^{\nu(t+s)} ds \leq C<\infty, \quad t\to \infty .
  \end{equation*}
Together with the boundedness of $I'$ and (\ref{0020}), this proves that
\begin{equation*}
I_{12}<C<\infty,\quad |x|\leq Bt^\gamma, ~ ~ t\to\infty.
 \end{equation*}

{\it 2) Estimate on $I_{11}$}.
The two balls $P_1,P_2$ do not intersect each other when  $ |x|>2t^{\frac{1}{\alpha}}$ (and $0\leq s\leq t$). Hence we may assume that $ |x|<2t^{\frac{1}{\alpha}}$. Then (\ref{t2ii}) implies that
$|m_1^2 (t,x)|\asymp   t^{-\frac{2d}{\alpha}} e^{2t\nu }, ~   t\to\infty $. We apply (\ref{t20}) to the factors under the summation sign in (\ref{0019}) and obtain that

\begin{equation*}
I_{11}\leq   \frac{C t^{\frac{2d}{\alpha}}}{ e^{2\nu t}}
  \int _{0}^{t}{
  \sum_{z\in D_{11} }{\frac{e^{\nu(t+s)}ds}{(t-s+1)^{\frac{d}{\alpha}}( s+1)^{\frac{2d}{\alpha}}}} }, \quad t\to\infty.
 \end{equation*}.

 We replace here the summation sign by the factor $\kappa$ introduced above and use the estimate
 \[
 \frac{ \kappa}{(t-s+1)^{\frac{d}{\alpha}}}\leq C\frac{ [(t-s)^{\frac{d}{\alpha}}+1]}{(t-s+1)^{\frac{d}{\alpha}}}\leq C_1.
 \]
This leads to
 \begin{equation}\label{0021}
  I_{11}\leq  \frac{C t^{\frac{2d}{\alpha}}}{ e^{\nu t}}
  \int_{0}^{t}
  {\frac{e^{\nu s}ds}{ (s+1)^{\frac{2d}{\alpha}}}}<C<\infty, \quad t\to\infty.
\end{equation}.

{\it 3) Estimate on $I_{22}$}. Estimate (\ref{t2ii}) can be applied to the factors under the summation sign in (\ref{0019}). Thus
\begin{equation} \label{1D22}
I_{22}\leq\frac{C}{ m_1^2(t,x)}
  \int_{0}^{t}{
  \sum_{z\in D_{22} } \frac{ (t-s)s^2}{|x-z|^{d+\alpha} |z|^{2(d+\alpha)}  }
 ~ e^{\nu(t+s)}}~ds.
 \end{equation}
 Consider first the case when $|x|<\frac{1}{2}t^{1/\alpha}$. Then the inequalities
 \begin{equation} \label{1111}
 |z|>s^\frac{1}{\alpha}, ~~~|x-z|>(t-s)^\frac{1}{\alpha},~ ~~0<s<t, ~~x,z\in \mathbb R^d,
 \end{equation}
  imply that $|z|\geq\beta>0$ if $t=1$. Using the homogeneity arguments, one can easily obtain that inequalities (\ref{1111}) with an arbitrary $t>0$ imply that
  $$|z|\geq\beta t^{1/\alpha},
  $$
 i.e., the latter estimate holds in $D_{22}.$ Now we can replace $|z|$ in (\ref{1D22}) by $\beta t^{1/\alpha}$ and use relation (\ref{t2ii}) for $m_1$. This leads to
\begin{equation*}
I_{22}\leq \frac{Ct^{2d/\alpha}}{ e^{2\nu t}t^{2(d+\alpha)/\alpha}}
  \int_{0}^{t}
  \sum_{z\in D_{22} } \frac{(t-s)s^2 e^{\nu(t+s)}ds}{|x-z|^{d+\alpha}   }, \quad |x|<\frac{1}{2}t^{1/\alpha}, ~~t\to\infty.
 \end{equation*}
Note that $|x-z|\geq 1$ when $z\in D_{22}$ and that $ \sum_{z\in D_{22} } \frac{1}{|x-z|^{d+\alpha}}<\sum_{z\in \mathbb Z^d\backslash\{x\} } \frac{1}{|x-z|^{d+\alpha}}$. The latter series converges and does not depend on $x$. Hence
\begin{equation} \label{1D22a}
I_{22}\leq \frac{Ct^{2}}{ e^{2\nu t}}
  \int_{0}^{t}
   (t-s)s^2 e^{\nu(t+s)}ds\leq C<\infty, \quad |x|<\frac{1}{2}t^{1/\alpha}, ~~t\to\infty.
 \end{equation}

Now let us estimate $I_{22}$ when $\frac{1}{2}t^{1/\alpha}\leq |x|\leq Bt^\gamma$. We split the region $D_{22}$ into two, namely,  $ D_{22}^{(1)}= D_{22 } \cap \{z:|z|> \frac{|x|}{2}\} $ and  $  D_{22}^{(2)}= D_{22} \cap \{z:|z|  \leq \frac{|x|}{2}  \} $. Then $ I_{22}\leq I_1+I_2$,
  where $I_1$ and $I_2$ are the right-hand side in (\ref{1D22}) with $D_{22} $  replaced by $D_{22}^{(1)} , D_{22}^{(2)}$ respectively.

 In order to evaluate $I_1$, we replace $|z|$ in (\ref{1D22}) by $|x|/2$, use convergence of the series $ \sum_{z\in D_{22} } \frac{1}{|x-z|^{d+\alpha}}$, and estimate (\ref{t2i}) for $m_1$. This leads to
\[
I_1\leq \frac{C|x|^{2( d+\alpha)}}{t^2 e^{2t\nu} }
  \int _{0}^{t} \frac{ (t-s)s^2}{|x|^{2( d+\alpha)}} e^{\nu(t+s)}ds=Ct^{-2}e^{-\nu t}
  \int _{0}^{t} s^2 (t-s) e^{\nu s}ds\leq C<\infty
\]
when $\frac{1}{2}t^{1/\alpha}\leq |x|\leq Bt^\gamma, ~~t\to\infty$.

Let us estimate $I_2$. The inequality $|z|\leq \frac{|x|}{2} $ implies that $|x-z|\geq \frac{|x|}{2} $. We use this in (\ref{1D22}) together with estimate (\ref{t2i}) for $m_1$. Then
we have
\begin{equation*}
I_2
 \leq \frac{C |x|^{(d+\alpha)}}{ t^2e^{\nu t}}
  \int_{0}^{t} (t-s) s^2 e^{\nu s}
 ( \sum_{z\in D_{22}^{(2)}} \frac{1}{|z|^{2(d+\alpha)}})ds, \quad  \frac{1}{2}t^{1/\alpha}\leq |x|\leq Bt^\gamma, ~~t\to\infty.
 \end{equation*}
Note that $|z|>s^\frac{1}{\alpha}$ in $D_{22}$ and $|z|\geq 1$ if $z\neq 0$. Hence the sum in the formula above does not exceed
$C(1+ s)^{-\frac{2\alpha+d}{\alpha}}$. We also can use there that $\frac{|x|^{(d+\alpha)}}{t^2}<C_1t^{d/\alpha}$ (due to $|x|\leq Bt^\gamma$). Hence
  \begin{equation*}
  I_2\leq  \frac{Ct^{d/\alpha}}{ e^{t\nu} }\int_0^t {(1+ s)^{-\frac{2\alpha+d}{\alpha}} (t-s) s^2  e^{\nu s}~ ds} , \quad  \frac{1}{2}t^{1/\alpha}\leq |x|\leq Bt^\gamma, ~~t\to\infty.
     \end{equation*}
 One can easily check that the right-hand side above is bounded. This together with the boundedness of $I_1$ and (\ref{1D22a}) implies that
 \[
 I_{22}\leq C<\infty, \quad |x|\leq Bt^\gamma, ~~t\to\infty.
 \]

{\it 4) Estimate on $I_{21}$}. From Lemma \ref{linf} it follows that
\[
I_{21}\leq\frac{C }{ m_1^2(t,x)}
  \int_{0}^{t}
  \sum_{z\in D_{21} } \frac{(t-s)}{|x-z|^{d+\alpha}  (s+1)^{\frac{2d}{\alpha}} }
 e^{\nu(t+s)}ds.
\]

First, let us estimate $I_{21}$ when $|x|\leq 2t ^{\frac{1}{\alpha }} $. Then $m_1$ can be estimated using (\ref{t2ii}). There is also $x$-independent constant $C<\infty$ such that $ \sum_{z\in D_{21} } \frac{1}{ |x-z|^{d+\alpha}   }<C$ (see details in the subsection on $I_{22}$). Thus
 \[
I_{21}\leq\frac{C t^{ \frac{2d}{\alpha}}}{ e^{2vt}}
  \int_{0}^{t}
  \frac{(t-s)}{ (s+1)^{\frac{2d}{\alpha}} }
 e^{\nu(t+s)}ds=\frac{C t^{ \frac{2d}{\alpha}}}{ e^{vt}}
  \int_{0}^{t}
  \frac{(t-s)}{ (s+1)^{\frac{2d}{\alpha}} }
 e^{\nu s}ds\leq C<\infty
\]
when $|x|\leq 2t ^{\frac{1}{\alpha }} , ~~t\to\infty.$

Consider now the case of $ 2t ^{\frac{1}{\alpha }}<|x|\leq Bt^\gamma.$ Then estimate  (\ref{t2i}) can be applied to $m_1$, and therefore
\begin{equation} \label{1D21}
I_{21}\leq \frac{C |x|^{2(d+\alpha)}}{ t^2e^{2\nu t}}
  \int_{0}^{t}
  \sum_{z\in D_{21} } \frac{(t-s) e^{\nu(t+s)}}{|x-z|^{d+\alpha}(s+1)^{\frac{2d}{\alpha}}
}ds.
 \end{equation}.

Since $|z|\leq s ^{\frac{1}{\alpha}} \leq t ^{\frac{1}{\alpha}}$ in $D_{21}$ and we assume that $ |x|>2t ^{\frac{1}{\alpha }}$, we have $|x-z|\geq |x|/2$ in the integrand above. Thus we can replace $|x-z|$ by $|x|/2$ in (\ref{1D21}). After that, the summation sign can be replaced by the number $\kappa_1$ of terms in the sum. Obviously, $\kappa_1\leq C(A_1^{d}+1)$, where $A_1$ is the radius of the first ball $P_1$, i.e., $ \kappa_1\leq C[s^{\frac{d}{\alpha}}+1],$ and
\begin{equation*}
I_{21}\leq \frac{C |x|^{d+\alpha}}{ t^2e^{\nu t}}
  \int_{0}^{t}
  \frac{(t-s) e^{\nu s}}{(s+1)^{\frac{d+1}{\alpha}}
}ds\leq \frac{C t^{d/\alpha}}{ e^{\nu t}}
  \int_{0}^{t}
  \frac{(t-s) e^{\nu s}}{(s+1)^{\frac{d+1}{\alpha}}
}ds\leq C<\infty
 \end{equation*}
when $2t ^{\frac{1}{\alpha }}<|x|\leq Bt^\gamma, ~~t\to\infty.$ Thus $I_{21}$ is bounded when $|x|\leq Bt^\gamma, ~~t\to\infty.$ Together with the boundedness of all other $I_{ij}$, this completes the proof of the theorem.

\qed

\bibliographystyle{amsalpha}
\bibliography{bibliography.bib}

\end{document}